\documentclass[12pt]{article}
\usepackage{amstex}
\usepackage{amssymb}
\usepackage{theorem}
\usepackage{a4}
\usepackage{pb-diagram}
\usepackage{lamsarrow}
\usepackage{pb-lams}

\def\omi#1{\buildrel { \buildrel{#1}\over{\vee} } \over .}

\def\bbbone{{ \mathchoice {1\mskip-4mu\mathrm{l} } {1\mskip-4mu\mathrm{l}}
{1\mskip-4.5mu\mathrm{l} } {1\mskip-5mu\mathrm{l}} }}

\def\gC{{\mathbb C}}

\def\gN{{\mathbb N}}

\def\cA{{\cal A}}

\def\cL{{\cal L}}

\def\cX{{\cal X}}
\def\cY{{\cal Y}}
\def\cZ{{\cal Z}}

\def\tA{{\mathfrak A}}

\def\fA{{\mathfrak A}}

\def\findem{~\hfill$\square$}

\def\demo{\noindent {\sl Proof: }}

\def\der{{\mathrm{Der}}}
\def\out{{\mathrm{Out}}}
\def\Int{{\mathrm{Int}}}
\def\End{{\mathrm{End}}}

\def\Id{{\mathrm{Id}}}

\def\tr{{\mathrm{Tr}}}

\def\exter{{\textstyle\bigwedge}}

\def\loc{{\mathrm{loc}}}

{\theorembodyfont{\sl}
\newtheorem{theo}{Theorem}
\newtheorem{prop}[theo]{Proposition}
\newtheorem{cor}[theo]{Corollary}

\newtheorem{lem}[theo]{Lemma}}
{\theorembodyfont{\rm}

}

\begin{document}
\baselineskip=0.7cm

\vspace*{3cm}
\begin{center} 
{\large\bf $SU(n)$-CONNECTIONS AND\\
\bigskip
NONCOMMUTATIVE DIFFERENTIAL GEOMETRY} 
\end{center}
\vspace{1cm}

\begin{center}
Michel DUBOIS-VIOLETTE and Thierry MASSON\\
\vspace{0.3cm}
{\small Laboratoire de Physique Th\'eorique et Hautes
Energies\footnote{Laboratoire associ\'e au Centre National de la
Recherche Scientifique - URA D0063}\\
Universit\'e Paris XI, B\^atiment 211\\
91 405 Orsay Cedex, France\\
e-mail: flad@@qcd.th.u-psud.fr and masson@@qcd.th.u-psud.fr}
\end{center}
\vspace{1cm}

\vspace{1cm}
\begin{abstract}
We study the noncommutative differential geometry of the algebra of
endomorphisms of any $SU(n)$-vector bundle. We show that ordinary
connections on such $SU(n)$-vector bundle can be interpreted in a
natural way as a noncommutative $1$-form on this algebra for the
differential calculus based on derivations. We interpret the Lie algebra
of derivations of the algebra of endomorphisms as a Lie algebroid. Then
we look at noncommutative connections as generalizations of these usual
connections.  \end{abstract}

\vfill \noindent L.P.T.H.E.-ORSAY 96/100\\ \vspace{0.5cm}

\newpage

The noncommutative differential geometry of the algebra of matrix valued
functions on a manifold has been studied in \cite{DVKM2}. There it was
pointed out that noncommutative connections are good candidate to unify
at the classical level ordinary gauge fields and Higgs fields in a
unique object. This idea has been widely used in a variety of contexts.
In this paper, as a generalization of \cite{DVKM2}, we consider the
noncommutative differential calculus based on derivations for the
algebra of endomorphisms of a $SU(n)$-vector bundle.

In Section~\ref{dbndc}, we recall the definition of the differential
calculus based on derivation and the main results for the algebra of
matrix valued functions. In Section~\ref{taevb}, we consider the algebra
of endomorphisms of a $SU(n)$-vector bundle. We show that the affine
space of $SU(n)$-connections is an affine subspace of the space of
noncommutative $1$-form of this algebra. In Section~\ref{derla}, the Lie
algebra of derivations of this algebra is interpreted as a Lie
algebroid. In Section~\ref{ncc}, noncommutative connections are
considered. It is shown that they incorporate naturally Higgs fields as
in the case of matrix valued functions.

\section{Derivations based noncommutative differential calculi}
\label{dbndc}

The differential calculus based on derivations is a natural
generalization the differential algebra of differential forms on a
manifold (\cite{MDV1, MDV2, DVM1} and references therein).

\subsection{Construction}

Let $\fA$ be an associative algebra with unit. Then one has the well
known results: \begin{lem} The vector space $\der(\fA)$ of derivations
of $\fA$ is a Lie algebra and a module over the center $\cZ(\fA)$ of
$\fA$. The center $\cZ(\fA)$ is stable by $\der(\fA)$. The vector space
of inner derivations $\Int(\fA)$ is a Lie ideal and a
$\cZ(\fA)$-submodule.  \end{lem}

The quotient $\der(\fA)/\Int(\fA)$ will be denoted by $\out(\fA)$. This
is then a Lie algebra and a module over $\cZ(\fA)$.

One can consider the complex $\underline{\Omega}_\der(\fA)$ of
$\cZ(\fA)$-multilinear antisymmetric maps from $\der(\fA)$ to $\fA$. It
is naturally a $\gN$-graded algebra. One defines a differential
$\hat{d}$ (of degree $1$) on this graded algebra by setting, for any
derivations $X_1, \dots, X_{n+1}$ and any $\omega \in
\underline{\Omega}^n_\der(\fA)$ \begin{eqnarray*} \hat{d}\omega(X_1,
\dots , X_{n+1}) &=& \sum_{i=1}^{n+1} (-1)^{i+1} X_i \omega( X_1, \dots
\omi{i} \dots, X_{n+1}) \nopagebreak\nonumber\\ \nopagebreak & & +
\sum_{1\leq i < j \leq n+1} (-1)^{i+j} \omega( [X_i, X_j], \dots \omi{i}
\dots \omi{j} \dots , X_{n+1}) \end{eqnarray*} The graded differential
algebra $(\underline{\Omega}_\der(\fA), \hat{d})$ is the first
noncommutative differential calculus we define on $\fA$.

The second one is the smallest differential subalgebra of
$\underline{\Omega}_\der(\fA)$ generated by $\fA$. We denote it by
$\Omega_\der(\fA)$. Any element of $\Omega^n_\der(\fA)$ is a sum of
terms $a_0 \hat{d} a_1 \dots \hat{d} a_n$ for $a_i \in \fA$.

\bigskip These definitions are generalizations of the usual differential
calculus on a manifold in the sense that when $\fA$ is the algebra
$C^\infty(M)$ of smooth functions on a finite dimensional regular smooth
manifold $M$, these differential algebras coincide with the graded
differential algebra of differential forms on $M$.

\subsection{The matrix algebra}

The noncommutative geometry of the matrix algebra for the differential
calculus based on derivations has been described in \cite{MDV2} and
\cite{DVKM1}. We briefly summarize here the results.

The algebra $M_n(\gC)$ has only inner derivations. The Lie algebra
$\der(M_n(\gC)) = \Int(M_n(\gC))$ can be identified with the Lie algebra
$sl(n,\gC)$. The two previously defined differential calculi are the
same and one has \[ \Omega_\der(M_n(\gC)) = M_n(\gC) \otimes \exter
sl(n,\gC)^\ast \] where $sl(n,\gC)^\ast$ is the dual of $sl(n,\gC)$. We
denote by $d'$ the differential on this complex.

There exists a particular $1$-form $\theta$ defined by \[
i\theta(ad_\gamma) = \gamma - {1\over n} \tr (\gamma)\bbbone\] for any
$\gamma \in M_n(\gC)$. This $1$-form satisfies \[ d' i\theta -
(i\theta)^2 = 0 \] and for any $\gamma \in M_n(\gC) =
\Omega^0_\der(M_n(\gC))$, one has $d' \gamma = [i\theta, \gamma]$.

\subsection{The matrix valued functions algebra}

The derivations based differential calculus for the algebra $\fA =
C^\infty(M)\otimes M_n(\gC)$ for a manifold $M$ has been studied in
\cite{DVKM2}. The main results are the following.

The center of the algebra $\fA$ is the algebra $C^\infty(M)$ of smooth
complex valued functions on $M$. The Lie algebra of derivations
$\der(\fA)$ split canonically as a $C^\infty(M)$-module into \[
\der(\fA) = [\der(C^\infty(M))\otimes \bbbone ] \oplus [ C^\infty(M)
\otimes \der(M_n(\gC)) ] \] where $\der(C^\infty(M))=\Gamma(TM)$ is the
ordinary Lie algebra of vector fields on $M$ and $\der(M_n(\gC)) =
sl(n,\gC)$.

This result implies the canonical decomposition \[ \Omega_\der(\fA) =
\Omega(M) \otimes \Omega_\der(M_n(\gC)) \] where $\Omega(M)$ is the
graded differential algebra of differential forms on $M$, with
differential $d$. The differential $\hat{d}$ on $\Omega_\der(\fA)$ is
the sum $\hat{d} = d + d'$. In particular, restricted to $\fA$, $\hat{d}
= d + ad_{i\theta}$.

The $1$-form $\theta$ is well defined in $\Omega^1_\der(\fA)$ if we
extend it on $\der(\fA)$ by zero on the $\Gamma(TM)$ terms. $\theta$ is
real in the sense that it is real-valued on real (i.e. hermitian)
derivations.

\section{The algebra of endomorphisms of a vector bundle} \label{taevb}

In the following, by a $SU(n)$-vector bundle $E$ we mean a hermitian
vector bundle of rank $n$ such that $\exter^n E$ is trivial (i.e.
trivializable with a given trivialization).

Let $E$ be a $SU(n)$-vector bundle over regular finite dimensional
smooth (i.e. paracompact, etc...) manifold $M$, and let us denote by
$\End(E)$ the fiber bundle of endomorphisms of $E$. The sections of this
fiber bundle in matrix algebras form a unital algebra, which we denote
by $\fA$. The hermitian structure gives a natural involution on this
algebra, which we denote by $S \mapsto S^\ast$. The center of this
algebra is exactly $C^\infty(M)$ (smooth functions on $M$ with values in
$\gC$), identifying $f\in C^\infty(M)$ with $f\bbbone \in \fA$.  The
trace map, defined on each fiber of $\End(E)$, gives a natural map \[
\tr : \fA \rightarrow C^\infty(M) \] Similarly, the determinant defines
a natural map \[ \det : \fA \rightarrow C^\infty(M) \]

By restriction to the center, there is a natural map \[ \rho : \der(\fA)
\rightarrow \der(C^\infty(M)) = \Gamma(TM) \] This map is the quotient
map in the short exact sequence of Lie algebras and
$C^\infty(M)$-modules \begin{equation} \label{sesder} 0 \rightarrow
\Int(\fA) \rightarrow \der(\fA) {\buildrel{\rho}\over \longrightarrow}
\out(\fA) = \Gamma(TM) \rightarrow 0 \end{equation} For any derivation
$\cX \in \der(\fA)$, let us denote by $X\in \Gamma(TM)$ the associated
vector field on $M$. Notice now that the $1$-form $i\theta$ defined in
the previous section is well defined here on $\Int(\fA)$: \[
i\theta(ad_\gamma) = \gamma - {1\over n}\tr(\gamma) \bbbone \] for any
$\gamma\in \fA$. Any inner derivative $ad_\gamma$ will be taken such
that $\gamma$ is traceless. It is then a section of the fiber bundle of
traceless endomorphisms of $E$. We denote by $\fA_0$ the set of such
sections. The Lie algebra $\Int(\fA)\simeq \fA_0$ operates in the sense
of H. Cartan on $\Omega_\der(\fA)$ \cite{CAR, MDV1}. The horizontal
elements for this operation can be considered as differential forms on
$M$ with values in $\End(E)$. As will be seen below, the curvature of a
connection on $E$ can be interpreted in this way. The basic forms are
ordinary differential forms on $M$. In the following, horizontality will
refer to this operation.

The subalgebra $\cZ(\fA)$ of $\fA$ can be considered as a quotient
manifold algebra in the sense of \cite{MAS}. In this interpretation, the
algebra $\fA$ looks like a principal bundle and $\cZ(\fA)$ as the
algebra of functions on the base space. The Lie algebra of the gauge
group, defined as the Lie algebra of derivations which are zero on
$\cZ(\fA)$, is then $\Int(\fA)$. As we will see later, this point of
view will be confirmed when we will consider ordinary connections on
$E$.

A derivation $\cX \in \der(\fA)$ will be said real if $(\cX S)^\ast =
\cX (S^\ast)$ for any $S \in \fA$. By duality, one defines hermitian and
antihermitian noncommutative differential forms.

\medskip Let us now look at the differential calculi based on
derivations.

\begin{prop} The two differential calculi $\Omega_\der(\fA)$ and
$\underline{\Omega}_\der(\fA)$ coincide.  \end{prop}

\medskip Indeed, by using a finite open covering of $M$ with open sets
which trivialize $E$, and results on smooth matrix valued functions, one
can show by standard arguments that any element of
$\underline{\Omega}^p_\der(\fA)$ can be written as a finite sum of terms
$a_0 \hat{d} a_1 \dots \hat{d} a_p$ with $a_i \in \fA$. It is then an
element of $\Omega^p_\der(\fA)$.

\medskip We will denote by $\hat{d}$ the differential on
$\Omega_\der(\fA) = \underline{\Omega}_\der(\fA)$.

\medskip In the trivial case ($E = M \times \gC^n$), the algebra $\fA$
is exactly $C^\infty(M) \otimes M_n(\gC)$, and we are back to the
situation of the previous section. In particular, the short exact
sequence (\ref{sesder}) splits canonically as $C^\infty(M)$-modules.

In the general case, one has the following situation.  \begin{prop}
\label{alpha} Let $\nabla^E$ be any connection on $E$. Then there exists
a noncommutative $1$-form $\alpha$ in $\Omega^1_\der(\fA)$ such that any
derivation $\cX \in \der(\fA)$ can be decomposed as \begin{equation}
\label{decder} \cX = \nabla_X - ad_{\alpha(\cX)} \end{equation} where
$\nabla$ is the naturally associated connection to $\nabla^E$ on the
fiber bundle $\End(E)$.\\ In other words, restricted to $\fA$, the
noncommutative differential $\hat{d}$ can be written with obvious
notation \[ \hat{d} = \nabla - ad_\alpha \] \end{prop}

We recall that $\nabla$ is the tensor product of the connections
$\nabla^E$ on $E$ and $\nabla^{E^\ast}$ on the dual vector bundle
$E^\ast$ of $E$ where $\nabla^{E^\ast}$ satisfies $X \langle \epsilon ,
e \rangle = \langle \nabla^{E^\ast}_X \epsilon , e \rangle + \langle
\epsilon , \nabla^E_X e \rangle$ for any sections $\epsilon$ of $E^\ast$
and $e$ of $E$.

\medskip \demo First, notice that for a fixed $X\in \Gamma(TM)$,
$\nabla_X : \fA \rightarrow \fA$ is naturally a derivation. Now, the
image in $\Gamma(TM)$ of the derivation $\cX - \nabla_X$ is zero, as can
be seen by applying it on any function in $C^\infty(M)$. So the
derivation $\cX - \nabla_X$ has to be in the kernel of the quotient map.
It is then inner. Define $\alpha$ by \[ \alpha(\cX) = - i\theta \left(
\cX - \nabla_X \right) \] This is obviously a noncommutative $1$-form on
$\fA$ which satisfies by definition $ad_{\alpha(\cX)} = \nabla_X -
\cX$.\findem

\bigskip Notice that this noncommutative $1$-form takes its values in
the traceless elements of $\fA$. In fact, one can see $\alpha$ as an
extension of $-i\theta$ to all derivations. Indeed, one has obviously
$\alpha(ad_\gamma) = -\gamma$. Recall the convention that $\tr(\gamma) =
0$.

This proposition gives us a splitting of the short exact
sequence~(\ref{sesder}) of $C^\infty(M)$-modules. This splitting is not
canonical and is only defined through a choice of a connection on $E$,
by the $C^\infty(M)$-linear map $X\mapsto \nabla_X$ from $\Gamma(TM)$
into $\der(\fA)$. This has to be compared with the usual (commutative)
situation where one can interpret a connection as a map from vector
fields on $M$ into vector fields on a principal bundle over $M$. In our
situation, the algebra $\fA$ plays a similar role to this principal
bundle.

\medskip Let us now look at this noncommutative $1$-form $\alpha$ in
local expressions.  On any open subset $U\subset M$ of a trivialization
of $\End(E)$, the algebra $\fA$ looks like $C^\infty(U)\otimes
M_n(\gC)$. Let $S\in \fA$ be a section. Restricted to $U$, this section
can be regarded as a matrix valued function $s : U \rightarrow
M_n(\gC)$. Let $U'\subset M$ be a second open subset of a trivialization
of $\End(E)$. The restriction to $U'$ of the section $S$ is an other
matrix valued function $s' : U' \rightarrow M_n(\gC)$ which is related
to $s$ by $s' = g^{-1} s g$ for a transition function $g : U\cap U'
\rightarrow SU(n)$.

A derivation $\cX$ of the algebra $\fA$ can also be trivialized. Because
the algebra looks locally as $C^\infty(U)\otimes M_n(\gC)$, this
trivialization has the form $\cX_\loc = X + ad_\gamma$ for $X$ an
ordinary vector fields on $U$ and $\gamma$ a traceless matrix valued
function on $U$. $X$ is nothing but the restriction to $U$ of the image
of $\cX$ in $\Gamma(TM)$. On $U'$ as above, the trivialization of $\cX$
has the form $\cX'_\loc = X + ad_{\gamma'}$. It is easy to show that on
$U\cap U'$, one then has \[ \gamma' = g^{-1}\gamma g + g^{-1} (Xg) \]
Because the structure group is $SU(n)$, the last term is traceless.

\begin{cor} The local expression of the noncommutative $1$-form $\alpha$
is \[ \alpha_\loc(\cX_\loc) = A(X) - \gamma \] where $A$ is the local
expression of the connection $1$-form of $\nabla^E$ (with values in the
Lie algebra $su(n)$).  \end{cor}

One can verify directly that those local expressions can be joined
together into the global objet $\alpha$ because the inhomogeneous term
of the transition relations of the $A$ fields cancels exactly the
inhomogeneous term of the transition relations of the $\gamma$ fields.
Notice the importance to work on a $SU(n)$-vector bundle. The traceless
of the connection $1$-form is essential here. Strictly speaking,
$SL(n,\gC)$-connections would be sufficient.

\begin{cor} The canonical mapping $\nabla^E \mapsto \alpha$ is an
isomorphism of affine spaces from the affine space of
$SU(n)$-connections on $E$ onto the affine space of traceless
antihermitian noncommutative $1$-forms on $\fA$ satisfying
$\alpha(ad_\gamma) = -\gamma$.  \end{cor}

\medskip Indeed, a noncommutative $1$-form $\alpha$ is the image of a
$SU(n)$-connection if and only if $\alpha(ad_\gamma) = -\gamma$ for any
$\gamma \in \fA_0$, $\tr \alpha(\cX) = 0$, and $\alpha(\cX)^\ast +
\alpha(\cX) = 0$ for any real $\cX\in \der(\fA)$. As we will see below,
this last condition can be interpreted as a compatibility condition with
a hermitian form. Now, one can work on the noncommutative $1$-form
$\alpha$ instead of $\nabla^E$.

\begin{prop} The Lie algebra of real derivations on $\fA$ acts on the
space of $SU(n)$-connections through the Lie derivative defined on
$\Omega_\der(\fA)$.\\ Restricted to inner real derivations, the Lie
derivative corresponds to infinitesimal gauge transformations on
connections.  \end{prop}

\demo It is easy to see that for any real derivation $\cX \in
\der(\fA)$, if $\alpha$ is the image of a $SU(n)$-connection on $E$,
then $\cL_\cX \alpha = (i_\cX \hat{d} + \hat{d}i_\cX)\alpha$ satisfies
\[ (\cL_\cX \alpha)(ad_\gamma) = 0 \] and \[ (\cL_\cX \alpha)(\cY)^\ast
+ (\cL_\cX \alpha)(\cY) = 0 \] for any $\gamma \in \fA_0$ and real $\cY
\in \der(\fA)$. Then $\alpha + \cL_\cX \alpha$ is the image of a
$SU(n)$-connection on $E$.  Now, for $\cX= ad_\xi$ with $\xi^\ast + \xi
= 0$, one has \[ \cL_{ad_\xi} \alpha = -\hat{d} \xi - [\alpha, \xi] \]
This is exactly an infinitesimal gauge transformation on $\nabla^E$ if
$\xi$ is considered as an element of the Lie algebra of the group of
gauge transformations of $E$.\findem

\medskip Notice that in this proof, we have used the fact that the real
inner derivations $ad_\xi$ are exactly in correspondence with the
elements $\xi$ of the Lie algebra of the gauge group of $E$.

\bigskip In ordinary differential geometry, the connection $\nabla^E$ on
$E$ is related to a covariant vector valued $1$-form on an associated
principal bundle (the bundle of $SU(n)$-frames of $E$). Here, as we have
seen above, the algebra $\fA$ replaces in a certain sense this principal
bundle. The connection $\nabla^E$ is then associated to a noncommutative
$1$-form on $\fA$. However, one has to notice two important facts.
First, there is no need to a covariance restriction. Instead, one fixes
the value of the $1$-form on inner derivations. Second, there is a
natural way to generalize this structure taking any noncommutative
$1$-form on $\fA$. In Section~\ref{ncc} we shall see that any
noncommutative $1$-form represents a noncommutative connection and gives
rise to a Yang-Mills-Higgs type multiplet.

\begin{prop} \label{courbure} Let $\alpha$ be the noncommutative
$1$-form associated to the connection $\nabla^E$ on $E$. Let us denote
by $R^E$ the curvature of $\nabla^E$. Then one has \[ R^E(X,Y) =
\hat{d}\alpha(\cX, \cY) + [\alpha(\cX), \alpha(\cY) ] \] for any $\cX,
\cY \in \der(\fA)$, $X,Y$ being their images in $\Gamma(TM)$.  In
particular, the expression $\hat{d}\alpha + \alpha^2$ is a horizontal
element of $\Omega^2_\der(\fA)$.  \end{prop}

\demo The direct computation of $\hat{d}\alpha(\cX, \cY) + [\alpha(\cX),
\alpha(\cY) ]$ in a local trivialization gives immediately $dA(X,Y) +
[A(X), A(Y)]$.\findem

\medskip This proposition gives us an interpretation of the curvature
$R^E$ in terms of a horizontal noncommutative $2$-form on $\fA$.

\bigskip The curvature of the associated connection $\nabla$ on
$\End(E)$ is $R = ad_{R^E}$. One verifies that the decomposition
$\nabla_X = \cX + ad_{\alpha(\cX)}$ yields directly to \[ R(X,Y) =
ad_{\hat{d}\alpha(\cX, \cY) + [\alpha(\cX), \alpha(\cY) ] } \]

The curvature of $\nabla^E$ can be interpreted as the obstruction on the
morphism of modules from vector fields on $M$ into vector fields on the
associated principal bundle, to be a morphism of Lie algebras.  The
above formula can be interpreted in a similar way: $R$ is the
obstruction on the application $X \mapsto \nabla_X$ to be a morphism of
Lie algebras from $\Gamma(TM)$ to $\der(\fA)$.

\section{$\der(\fA)$ as a Lie algebroid} \label{derla}

There is a deep connection between the above discussion and the notion
of Lie algebroid. Let us recall first the definition of Lie algebroids
\cite{MAC}.

A Lie algebroid is a vector bundle $L$ over a smooth manifold $M$ with a
structure of Lie algebra on its smooth sections and a vector bundle
morphism $\rho : L \rightarrow TM$, called the anchor, such that \[ [
\rho(a) , \rho(b) ] = \rho([a,b]) \] and \[ [a, fb] = f[a,b] + (\rho(a)
f) b \] where $a$ and $b$ are sections of $L$, $[, ]$ is the Lie bracket
on sections of $L$ (or on vector fields on $M$ in the first formulae)
and $f$ is a function on $M$ (the anchor being naturally defined on
sections of $L$).

One important example of Lie algebroid is the Atiyah algebra $\cA(E)$
associated to any vector bundle $E$ over $M$. It is defined as the Lie
algebra of first-order differential operators on $E$ with symbol $\Id_E
\otimes X$ for $X$ a vector field on $M$, the anchor being the symbol
map $\sigma : \cA(E) \rightarrow \Gamma(TM)$. One has a natural short
exact sequence \[ 0 \rightarrow \Gamma( \End(E)) \rightarrow \cA(E)
{\buildrel{\sigma}\over\longrightarrow} \Gamma(TM) \rightarrow 0 \] of
sections of vector bundles.

Associated to any Lie algebroid $L$ there is a differential calculus on
the sections of $\exter L^\ast$ with differential $d_L$ given by
\begin{eqnarray} \label{dl} d_L \omega(\ell_1, \dots , \ell_{n+1}) &=&
\sum_{i=1}^{n+1} (-1)^{i+1} \rho(\ell_i) \omega( \ell_1, \dots \omi{i}
\dots, \ell_{n+1}) \nopagebreak\nonumber\\ \nopagebreak & & +
\sum_{1\leq i < j \leq n+1} (-1)^{i+j} \omega( [\ell_i, \ell_j], \dots
\omi{i} \dots \omi{j} \dots , \ell_{n+1}) \end{eqnarray} for any $\omega
\in \exter^n L^\ast$ and $\ell_1, \dots, \ell_{n+1} \in L$.

With previous notations, a {\sl $L$-connection on $E$} is an anchor
preserving map \[ D : L \rightarrow \cA(E) \] The curvature of such a
connection is defined as the obstruction to be a Lie algebras morphism.

\bigskip Consider now $E$ as in the previous section. Then $\der(\fA)$
is in a natural way the sections of a Lie algebroid, with anchor given
by the restriction $\rho$ to $C^\infty(M) = \cZ(\fA)$ of derivations of
$\fA$ (see the sequence (\ref{sesder})). It is easy to see that
$\underline{\Omega}_{\der}(\fA)$ is just $\fA \otimes_{\cZ(\fA)}
\exter_{\cZ(\fA)} \der(\fA)^\ast$ and the differential associated to the
Lie algebroid $\der(\fA)$ defined on $\exter_{\cZ(\fA)} \der(\fA)^\ast$
as above (\ref{dl}) is just the restriction of $\hat{d}$.

There is a natural anchor preserving map $\cA(E) \rightarrow \der(\fA)$
which associates to any $T \in \cA(E)$, the derivation $S \mapsto [T,
S]$ for any $S\in \fA$ where the commutator is that of operators on $E$.
Locally, it is also given by $\Id_E \otimes X + A \mapsto X + ad_A$ with
obvious notations. This map, together with previously defined ones,
gives us the following commutative exact diagram of
$C^\infty(M)$-modules and Lie algebras:

\[ \begin{diagram} \node[2]{0} \arrow{s} \node{0} \arrow{s} \\ \node{0}
\arrow{e} \node{C^\infty(M)} \arrow{s} \arrow{e,=} \node{C^\infty(M)}
\arrow{s} \arrow{e} \node{0} \arrow{s} \\ \node{0} \arrow{e} \node{\fA}
\arrow{s,l}{ad} \arrow{e,t}{i} \node{\cA(E)} \arrow{e,t}{\sigma}
\arrow{s} \node{\Gamma(TM)} \arrow{s,=} \arrow{e} \node{0} \\ \node{0}
\arrow{e} \node{ \Int(\fA) } \arrow{s} \arrow{e} \node{\der(\fA)}
\arrow{e,t}{\rho} \arrow{s} \node{\Gamma(TM)} \arrow{s} \arrow{e}
\node{0} \\ \node[2]{0} \node{0} \node{0} \end{diagram} \]

In this diagram, $i\theta : \Int(\fA) \rightarrow \fA$ gives a splitting
as $C^\infty(M)$-modules and Lie algebras of the first column. Any
ordinary connection on $E$ splits the short exact sequence of the second
line by $X \mapsto \nabla^E_X$ and of the third line by $X \mapsto
\nabla_X$. Any noncommutative $1$-form is a map from $\der(\fA)$ to
$\fA$. So, the difference of two $\der(\fA)$-connections in the sense of
Lie algebroids is a noncommutative $1$-form on $\fA$.

\begin{prop} For any given $SU(n)$-connection $\nabla^E$ on $E$ and any
$\cX \in \der(\fA)$, define ${\buildrel{o}\over D} : \der(\fA)
\rightarrow \cA(E)$ by \begin{equation} \label{defdo} {\buildrel{o}\over
D}(\cX) = \nabla^E_X - \alpha(\cX) \end{equation} Then one has the
following:\\ (i) $\buildrel{o}\over D$ is independent of the choice of
$\nabla^E$,\\ (ii) ${\buildrel{o}\over D}$ is a splitting of the short
exact sequence \begin{equation} \label{sesae} 0 \rightarrow C^\infty(M)
\rightarrow \cA(E) \rightarrow \der(\fA) \rightarrow 0 \end{equation}
considered as an exact sequence of $C^\infty(M)$-modules as well as an
exact sequence of Lie algebras,\\ (iii) ${\buildrel{o}\over D}$ induces
the splitting $i\theta$ of \begin{equation} \label{sesfa} 0 \rightarrow
C^\infty(M) \rightarrow \fA \rightarrow \Int(\fA) \rightarrow 0
\end{equation} \end{prop}

It follows in particular that $\buildrel{o}\over D$ is a
$\der(\fA)$-connection with vanishing curvature.

\demo For any given $SU(n)$-connection $\nabla^E$ on $E$, with
definition~(\ref{defdo}), for any $\cX\in \der(\fA)$ and any $S\in \fA$,
one has $[{\buildrel{o}\over D}(\cX) , S] = \nabla_X S - [\alpha(\cX),
S] = \cX S$. So ${\buildrel{o}\over D}$ is a splitting of (\ref{sesae})
as $C^\infty(M)$-modules. For any $\cX, \cY \in \der(\fA)$, one can show
by direct computation, using Proposition~\ref{courbure}, that $[
{\buildrel{o}\over D}(\cX) , {\buildrel{o}\over D}(\cY) ] -
{\buildrel{o}\over D}([\cX, \cY]) = 0$. So ${\buildrel{o}\over D}$ is a
splitting of (\ref{sesae}) as Lie algebras. Restricted to inner
derivations, ${\buildrel{o}\over D}$ coincides with $i\theta$. Then
${\buildrel{o}\over D}$ induces the splitting of (\ref{sesfa}).\\ If
${\nabla'}^E$ is an other $SU(n)$-connection on $E$, then ${\nabla'}_X^E
- \alpha'(\cX) - ( \nabla^E_X - \alpha(\cX) ) = {\nabla'}_X^E -
\nabla^E_X - i\theta( \nabla'_X - \nabla_X )$, which is zero because of
the relations between $\nabla^E$ and $\nabla$.\findem

Is is obvious that ${\buildrel{o}\over D}$ is anchor preserving, and is
then a $\der(\fA)$-connection. Because it is a splitting of
(\ref{sesae}) as Lie algebras, its curvature vanishes.

\medskip Notice that in a local trivialization of $E$,
${\buildrel{o}\over D}$ is given by $X + ad_\gamma \mapsto \Id_E \otimes
X + \gamma$ where as usual, $\gamma$ is traceless.

\medskip Any other splitting of (\ref{sesae}) as $C^\infty(M)$-modules
and Lie algebras is of the form $\cX \mapsto {\buildrel{o}\over D}(\cX)
+ \phi(X)$ with $\phi \in \Omega^1(M)$ and $d\phi = 0$. In particular,
it always induces the splitting $i\theta$. If one drops out the
splitting as Lie algebras, then any other splitting of (\ref{sesae}) as
$C^\infty(M)$-modules is of the form $\cX \mapsto {\buildrel{o}\over
D}(\cX) + \phi(\cX)$ with $\phi$ a noncommutative $1$-form on $\fA$ with
values in $C^\infty(M)$. The curvature of such a splitting is
$\hat{d}\phi$.  \medskip Then ${\buildrel{o}\over D}$ identifies
$\der(\fA)$ as a Lie subalgebroid of $\cA(E)$. This Lie subalgebroid is
the set of elements of $\cA(E)$ which preserve the volume on $E$. In
this identification, $\Int(\fA)$ is mapped to $\fA_0$.

\section{Noncommutative connections} \label{ncc}

Let us now turn to noncommutative connections for the differential
calculus $\Omega_\der(\fA)$ on the previously defined algebra $\fA$. We
consider here connections on right or left $\fA$-modules. These
noncommutative connections have been defined in \cite{CO1} and used in
many approaches to noncommutative geometry.

\subsection{The right module $\fA$}

The simplest example is the right module $\fA$ itself which is the free
right $\fA$-module of rank $1$. On this right module there is a natural
hermitian form defined by $\langle S, S' \rangle = S^\ast S' \in \fA$
for any $S,S'\in \fA$. In this case, a noncommutative connection is an
application \[ \widehat{\nabla}_\cX : \fA \rightarrow \fA \] such that
$\widehat{\nabla}_\cX(S S') = S \cX(S') + \widehat{\nabla}_\cX(S) S'$
and $\widehat{\nabla}_{f\cX} S = f\widehat{\nabla}_\cX S$ for any $S,S'
\in \fA$ and $f\in \cZ(\fA)$. The curvature of a noncommutative
connection is defined by $\hat{R}(\cX, \cY) S = [ \widehat{\nabla}_\cX,
\widehat{\nabla}_\cY ] S - \widehat{\nabla}_{[\cX, \cY]}S$ for any $S
\in \fA$ and $\cX, \cY \in \der(\fA)$, which is a right $\fA$-module
homomorphism.

Any noncommutative right-connection $\widehat{\nabla}$ on $\fA$ is
entirely given by $\widehat{\nabla}_\cX \bbbone = \omega(\cX)$ where
$\omega$ is a noncommutative $1$-form in $\Omega^1_\der(\fA)$. One can
then write \[ \widehat{\nabla}_\cX S = \cX S + \omega(\cX) S \] for any
$S\in \fA$. The expression $\widehat{\nabla}^0_\cX S = \cX S$ defines a
noncommutative connection on $\fA$ whose curvature is $0$. This is a
particular point in the affine space of noncommutative connections on
the right module $\fA$. It is easy to see that the curvature of
$\widehat{\nabla}$ is the left multiplication by the noncommutative
$2$-form \[ \hat{d}\omega (\cX, \cY) + [ \omega(\cX), \omega(\cY) ] \]

A connection is said to be compatible with the hermitian structure if \[
\cX \langle S, S' \rangle = \langle \widehat{\nabla}_\cX S, S' \rangle +
\langle S, \widehat{\nabla}_\cX S' \rangle \] for any $S,S' \in \fA$ and
real $\cX \in \der(\tA)$.  This compatibility condition is equivalent to
\[ \omega(\cX)^\ast + \omega(\cX) = 0 \] for any real $\cX \in
\der(\fA)$.

Any unitary element $U\in \fA$ with $\det(U) = \bbbone$ defines on $\fA$
a right module endomorphism which preserves the hermitian structure and
the $\det$ application by setting $S \mapsto US$. We denote by $SU(\fA)$
the group of such elements of $\fA$. This is exactly the gauge group of
the $SU(n)$-vector bundle $E$. We denote by $U(\fA)$ the group of
unitary elements of $\fA$. For any $U\in U(\fA)$, the gauge
transformation of a noncommutative connection $\widehat{\nabla}$ is
defined by $\widehat{\nabla}^U_\cX S = U^\ast \widehat{\nabla}_\cX
(US)$. The noncommutative $1$-form $\omega$ is then transformed as \[
\omega \mapsto U^\ast \omega U + U^\ast \hat{d} U \]

The next proposition says that any (commutative) connection on $E$
defines canonically a noncommutative connection on $\fA$.  \begin{prop}
For a fixed choice of a $SU(n)$-connection $\nabla^E$ on $E$, then of
$\alpha$ as in Proposition~\ref{alpha}, one defines a noncommutative
connection $\widehat{\nabla}^\alpha$ by setting \[
\widehat{\nabla}_\cX^\alpha S = \nabla_X S + S \alpha(\cX) = \cX S +
\alpha(\cX)S \] for any $\cX \in \der(\fA)$ and $S\in \fA$.\\ The
curvature of this connection is $\hat{R}^\alpha(\cX, \cY) = R^E(X,Y)$.\\
This noncommutative connection $\widehat{\nabla}^\alpha$ is compatible
with the hermitian structure on $\fA$.\\ A gauge transformation on
$\nabla^E$ induces a $SU(\fA)$-gauge transformation on
$\widehat{\nabla}^\alpha$.  \end{prop}

\demo This are just computations using the properties of $\alpha$ and
$\nabla$.\findem

\medskip These noncommutative connections are then other particular
points in the space of noncommutative connections. Any other
noncommutative connection $\widehat{\nabla}$ can be decomposed as \[
\widehat{\nabla}_\cX S = \widehat{\nabla}_\cX^\alpha S + \cA(\cX) S \]
with $\cA\in \Omega^1_\der(\fA)$. In this case, one has \[ \omega(\cX) =
\alpha(\cX) + \cA(\cX) \] Using (\ref{decder}), the supplementary term
$\cA$ can be decomposed itself as \[ \cA(\cX) = a(X) - B(\alpha(\cX)) \]
Then $B$ is an application \[ B : \fA_0 \simeq \Int(\fA) \rightarrow \fA
\] defined by $B(\gamma) = \cA(ad_\gamma)$ for any $\gamma\in \fA_0$.

Physically, the $B$ term has been studied in \cite{DVKM2} for the
particular algebra $C^\infty(M)\otimes M_n(\gC)$. It was claim there
that it can be interpreted as Higgs fields. In the present approach, one
can see that $B$ is independent of the choice of the connection
$\nabla^E$ on $E$, because it can be written $B(\gamma) =
\omega(ad_\gamma) + \gamma$. The field $B$ is then a well defined object
associated to $\widehat{\nabla}$ and corresponds to the ``purely
noncommutative'' part of the connection in the trivial situation
$C^\infty(M)\otimes M_n(\gC)$. In the general case, it is impossible to
canonically split the connection into a ``commutative'' part and a
``purely noncommutative'' one. One has to choose an ordinary connection
on $E$ to decompose it.

In a gauge transformation, one wishes to fix the connection of reference
$\alpha$. In this case, the transformations relations of $\cA$, $a$ and
$B$ are given by \begin{eqnarray*} \cA &\mapsto& U^\ast \cA U + U^\ast
(\nabla U) \\ a &\mapsto& U^\ast a U + U^\ast (\nabla U) \\ B &\mapsto&
U^\ast B U \end{eqnarray*} which are almost the ordinary gauge
transformations of the connection $1$-form for $\cA$ and $a$, but with
the ordinary differential replaced by $\nabla$.

\bigskip In terms of $\alpha$ and $\cA$, the curvature $\hat{R}$ of
$\widehat{\nabla}$ is given by \begin{eqnarray*} \hat{R}(\cX, \cY) &=&
\hat{R}^\alpha(\cX, \cY) \\ & & + \nabla_X \cA(\cY) - \nabla_Y \cA(\cX)
- \cA([\cX,\cY]) + [\cA(\cX), \cA(\cY)] \\ &=& R^E(X, Y) \\ & & +
(\hat{d}\cA)(\cX, \cY) + [ \cA(\cX), \cA(\cY) ] \\ & & + [\alpha(\cX),
\cA(\cY) ] - [ \alpha(\cY), \cA(\cX) ] \end{eqnarray*} With the previous
decomposition, one then has \begin{eqnarray*} \hat{R}(\cX, \cY) &=&
R^E(X, Y) \\ & & + \nabla_X a(Y) - \nabla_Y a(X) - a([X,Y]) + [a(X),
a(Y)] \\ & & - \nabla_X B(\alpha(\cY)) - [a(X), B(\alpha(\cY)) ] \\ & &
+ \nabla_Y B(\alpha(\cX)) + [a(Y), B(\alpha(\cX)) ] \\ & & + [
B(\alpha(\cX)), B(\alpha(\cY))] + B(\alpha([\cX, \cY])) \end{eqnarray*}
Notice that the term $R^E(X, Y) + \nabla_X a(Y) - \nabla_Y a(X) -
a([X,Y]) + [a(X), a(Y)]$ can be interpreted as the curvature of the
connection ${\nabla'}^E_X = \nabla^E_X + a(X)$ on $E$. Then the third
and the fourth lines are expressions of the type ${\nabla'}_X
B(\alpha(\cY))$.

\medskip In \cite{DVKM2}, for the trivial case $C^\infty(M)\otimes
M_n(\gC)$, it was shown that the minima of the action written with this
curvature is related to its horizontality. Then in the general case, let
us look at the horizontality of the curvature $\hat{R}$. Taking $\cX =
ad_\gamma$ and $\cY = ad_\eta$ one gets the condition \begin{equation}
\label{cond1} [ B(\gamma), B(\eta) ] - B([\gamma, \eta]) = 0
\end{equation} and taking any $\cY$ one gets, using the previous
relation, \begin{equation} \label{cond2} \nabla'_Y B(\gamma) - B(
\nabla_Y \gamma) = 0 \end{equation} which can be written \[ \nabla_Y B +
[a(Y), B] = 0 \] In case these two relations are satisfied, $\hat{R}$ is
horizontal and \[ \hat{R}(\cX, \cY) = {R'}^E(X,Y) - B( R^E(X,Y)) \] The
relations (\ref{cond1}) and (\ref{cond2}) are the generalizations of the
expressions found in \cite{DVKM2} for the trivial case.

Notice that this computation has been done using the decomposition of
the noncommutative $1$-form $\cA$ in $\omega = \alpha + \cA$. It is
possible to decompose directly the noncommutative $1$-form $\omega$ as
$\omega(\cX) = a_\omega(X) - B_\omega(\alpha(\cX))$ and then compute the
curvature. But the expressions of the horizontality conditions on
$a_\omega$ and $B_\omega$ found in this way are less suggestive than
(\ref{cond1}) and (\ref{cond2}).

\bigskip One way to naturally generalize the previous discussion is to
consider a right module of the type $\fA^m$ for $m\in \gN$. This leads
to similar results with much more structures, in particular in the Higgs
sector (see \cite{DVKM2} for the trivial case).

\subsection{The left module $\Gamma(E)$}

As a second canonical example, consider the space $\Gamma(E)$ of
sections of $E$ which is a natural left $\fA$-module. Then one has the
following result: \begin{prop} There is a one-to-one correspondence
between splittings of (\ref{sesae}) as $C^\infty(M)$-modules and
noncommutative left-connections on $\Gamma(E)$. The curvature of a
noncommutative left-connections on $\Gamma(E)$ corresponds to the
obstruction to be a Lie algebra splitting.  \end{prop}

\demo If $\widetilde{\nabla}$ is a noncommutative left-connection on
$\Gamma(E)$, then $\widetilde{\nabla}_\cX$ is a first order operator on
$E$. Using the relation $\widetilde{\nabla}_\cX (Se) = (\cX S)e +
S\widetilde{\nabla}_\cX e$ for any $e\in \Gamma(E)$ and $S \in \fA$, one
can show that $[\widetilde{\nabla}_\cX , S ] = \cX S$. So $\cX \mapsto
\widetilde{\nabla}_\cX$ is a splitting of (\ref{sesae}) as
$C^\infty(M)$-modules. If $D : \der(\fA) \rightarrow \cA(E)$ is a
splitting of (\ref{sesae}) as $C^\infty(M)$-modules then $D(\cX)(Se) =
[D(\cX), S] e + S D(\cX)e$ for any $S \in \fA$ and $e \in \Gamma(E)$.
The first term is $(\cX S) e$, and so $D$ satisfies the derivation rule
of a noncommutative left-connection on $\Gamma(E)$.\findem

\medskip Then the canonical splitting ${\buildrel{o}\over D}$ defined in
Section~\ref{derla} gives us a canonical noncommutative connection of
vanishing curvature \[ {\buildrel{o}\over D}(\cX) e = \nabla^E_X e -
\alpha(\cX) e \] for any $e\in \Gamma(E)$ and $\cX \in \der(\fA)$.

\bigskip In this case, it is also possible to consider generalizations
taking as left module a direct sum of $\Gamma(E)$.

\section{Conclusions}

Here we have analyzed the noncommutative differential geometry of the
algebra $\fA$ of sections of the endomorphism bundle of a $SU(n)$-vector
bundle, thereby generalizing several results on the algebra of matrix
valued functions. One advantage of doing this is to isolate what is
canonical and what depends on the choice of a connection. It is worth
noticing here that not every bundle in matrix algebras can be identified
as the endomorphism bundle of a vector bundle (there are well known
homological obstructions of that).

We have shown that connections on $SU(n)$-vector bundles identify with
the affine subspace of $\Omega_\der(\fA)$ of traceless antihermitian
elements satisfying $\alpha(ad_\gamma) = -\gamma$. This last condition
is the analog of the vertical projection property of the corresponding
connection form on the appropriate principal $SU(n)$-bundle. Notice
however that here no equivariant property is required. In fact, it is
apparent here that the noncommutative algebra $\fA$ can be used in many
respects like this principal bundle. Doing that, we also make a bridge
between the notion of Lie algebroids and the noncommutative differential
calculus based on derivations.  Furthermore, we have shown that the
noncommutative connections on right or left modules are natural
extensions of the usual connections.

Here we have restricted attention to the structure group $SU(n)$. But it
is clear that by forgetting hermitian properties, one can pass to $SL(n,
\gC)$, and by imposing reality, one can pass to $SO(n)$.

\end{document}